\documentclass[11pt,leqno]{amsart}

\usepackage[scaled=0.92]{helvet}	

\usepackage[margin=.9in]{geometry}                
\usepackage{graphicx,multicol,color}
\usepackage{amssymb}
\usepackage{pstricks,pst-plot}

\usepackage{mathptmx}

\newcommand{\Ad}{\operatorname{Ad}}
\newcommand{\C}{\mathbb{C}}
\newcommand{\comp}{\mathrel{\scriptstyle\circ}}
\newcommand{\g}{\mathfrak{g}}
\newcommand{\inv}{^{-1}}
\newcommand{\pt}{{\rm pt}}
\newcommand{\R}{\mathbb{R}}
\newcommand{\vol}{{\rm vol}}
\newcommand{\Zero}{{\rm Zero}}

\usepackage{amsmath}		

\begin{document}
\title{What Is $\ldots$ Equivariant Cohomology?}
\author{Loring W. Tu}

\begin{multicols}{2}
\thanks{Loring W. Tu is Professor of Mathematics at Tufts University 
in Medford, Massachusetts.  His email address is 
{\tt loring.tu@tufts.edu}.}
\maketitle

\section*{Introduction}

Many invariants in geometry and topology can be
computed as integrals.
For example, in classical differential geometry the Gauss--Bonnet theorem states that if $M$ is a compact, oriented surface in Euclidean $3$-space with Gaussian curvature $K$
and volume form $\vol$, then its Euler characteristic is
\[
\chi(M) = \frac{1}{2\pi} \int_M K\,\vol.
\]
On the other hand, if there is a continuous vector field $X$ with
isolated zeros on a compact, oriented manifold, the Hopf index
theorem in topology computes the Euler characteristic of the manifold
as the sum of the indices at the zeros of the vector field $X$.
Putting the two theorems together, one obtains
\begin{equation} \label{e:gbh}
\frac{1}{2\pi}\int_M K\, \vol = \sum_{p \in \Zero(X)} i_X(p),
\end{equation}
where $i_X(p)$ is the index of the vector field $X$ at the zero $p$.
This is an example of a localization formula, for it computes a global
integral in terms of local information at a finite set of points. 
More generally, one might ask what kind of integrals can be computed
as finite sums.
A natural context for studying this problem is the situation when
there is a group acting on the manifold with isolated fixed points.
In this case, one can try to relate an integral over the manifold to a sum over the fixed point set.

Rotating the unit sphere $S^2$ in $\R^3$ about the $z$-axis
is an example of an action of the circle $S^1$ on 
the sphere.
It has two fixed points, the north pole and the south pole.
This circle action generates a continuous vector field on the sphere,
and the zeros of the vector field are precisely the fixed points of
the action (see Figure~1). 

\begin{center}
{\psset{unit=.8}
\begin{pspicture}(-2,-2.6)(2,3)
\pscircle(0,0){2}
\parametricplot[linestyle=dashed]{0}{180}{t cos 2 mul t sin .6 mul}
\parametricplot{180}{360}{t cos 2 mul t sin .6 mul}
\parametricplot[linewidth=2.5pt,arrows=->]{230}{250}{t cos 2 mul t sin .6 mul}
\parametricplot[linewidth=2.5pt,arrows=->]{260}{280}{t cos 2 mul t sin .6 mul}
\parametricplot[linewidth=2.5pt,arrows=->]{290}{310}{t cos 2 mul t sin .6 mul}
\rput(0,.5){
\parametricplot[linewidth=.5pt]{180}{360}{t cos 1.936 mul t sin .581 mul}
\parametricplot[linewidth=2.5pt,arrows=->]{260}{280}{t cos 1.936 mul t sin .581 mul}
\parametricplot[linewidth=2.5pt,arrows=->]{230}{250}{t cos 1.936 mul t sin .581 mul}
\parametricplot[linewidth=2.5pt,arrows=->]{290}{310}{t cos 1.936 mul t sin .581 mul}
}
\rput(0,1){
\parametricplot[linewidth=.5pt]{180}{360}{t cos 1.732 mul t sin .520 mul}
\parametricplot[linewidth=2.5pt,arrows=->]{230}{250}{t cos 1.732 mul t sin .520 mul}
\parametricplot[linewidth=2.5pt,arrows=->]{260}{280}{t cos 1.732 mul t sin .520 mul}
\parametricplot[linewidth=2.5pt,arrows=->]{290}{310}{t cos 1.732 mul t sin .520 mul}
}
\rput(0,1.5){
\parametricplot[linewidth=.5pt]{180}{360}{t cos 1.323 mul t sin .397 mul}
\parametricplot[linewidth=2.5pt,arrows=->]{230}{250}{t cos 1.323 mul t sin .397 mul}
\parametricplot[linewidth=2.5pt,arrows=->]{260}{280}{t cos 1.323 mul t sin .397 mul}
\parametricplot[linewidth=2.5pt,arrows=->]{290}{310}{t cos 1.323 mul t sin .397 mul}
}
\rput(0,-.5){
\parametricplot[linewidth=.5pt]{180}{360}{t cos 1.936 mul t sin .581 mul}
\parametricplot[linewidth=2.5pt,arrows=->]{260}{280}{t cos 1.936 mul t sin .581 mul}
\parametricplot[linewidth=2.5pt,arrows=->]{230}{250}{t cos 1.936 mul t sin .581 mul}
\parametricplot[linewidth=2.5pt,arrows=->]{290}{310}{t cos 1.936 mul t sin .581 mul}
}
\rput(0,-1){
\parametricplot[linewidth=.5pt]{180}{360}{t cos 1.732 mul t sin .520 mul}
\parametricplot[linewidth=2.5pt,arrows=->]{230}{250}{t cos 1.732 mul t sin .520 mul}
\parametricplot[linewidth=2.5pt,arrows=->]{260}{280}{t cos 1.732 mul t sin .520 mul}
\parametricplot[linewidth=2.5pt,arrows=->]{290}{310}{t cos 1.732 mul t sin .520 mul}
}
\psdots(0,2)(0,-2)
\rput(0,-2.7){Figure 1. A circle action on a sphere.}
\end{pspicture}
}
\end{center}

\bigskip

Recall that the familiar theory of singular cohomology 
gives a functor from the category of topological spaces 
and continuous maps to the category of graded rings 
and their homomorphisms.
When the topological space has a group action, one would like a functor that reflects both the topology of the space and the action of the group.
Equivariant cohomology is one such functor.

The origin of equivariant cohomology is somewhat convoluted.
In 1959 Borel defined equivariant singular cohomology in
the topological category using a construction now called
the Borel construction.
Nine years earlier, in 1950,
in two influential articles on the
cohomology of a manifold $M$ acted on by
a compact, connected Lie group $G$,
Cartan constructed a differential complex
$\left(\Omega_G^*(M), d_G\right)$
 out of the differential forms 
on $M$ and the Lie algebra of $G$.
Although the term ``equivariant cohomology'' never occurs in
Cartan's papers, Cartan's complex turns out to
 compute the real equivariant singular cohomology of a
$G$-manifold (a manifold with an action of a
Lie group $G$), in much the same way that the
de Rham complex of smooth differential forms computes the real
singular cohomology of a manifold.
Without explicitly stating it, Cartan provided the key step in
a proof of the equivariant de Rham theorem,
before equivariant cohomology was even defined!
In fact, a special case of the Borel construction was already
present in Cartan's earlier article (Colloque de Topologie, C.B.R.M.,
Bruxelles, 1950, p.~62).
Elements of Cartan's complex are called \emph{equivariant differential
  forms} or \emph{equivariant forms}.
Let $S(\g^*)$ be the polynomial algebra on the Lie algebra $\g$ of $G$;
it is the algebra of all polynomials in linear forms on $\g$.
An equivariant form on a $G$-manifold $M$ is a differential form
$\omega$ on $M$ with values in the polynomial algebra $S(\g^*)$
satisfying the equivariance condition:
\[
\ell_g^* \omega = (\Ad g\inv) \comp \omega \quad \text{for all\ } g \in G,
\]
where $\ell_g^*$ is the pullback by left multiplication by $g$
and $\Ad$ is the adjoint representation.
An equivariant form $\omega$ is said to be \emph{closed} if it
satisfies $d_G \omega = 0$.

What makes equivariant cohomology particularly useful
in the computation of integrals is the equivariant
integration formula of Atiyah--Bott (1984) and 
Berline--Vergne (1982).
In case a torus acts on a compact, oriented manifold
with isolated fixed points,
this formula computes the integral of a closed equivariant form as a
finite sum over the fixed point set.
Although stated in terms of equivariant cohomology, 
the equivariant integration formula, also called the equivariant localization formula in the literature, can often be used to compute the
integrals of ordinary differential forms.
It opens up the possibility of machine computation of integrals on a
manifold. 

\section*{Equivariant Cohomology}

Suppose a topological group $G$ acts continuously on a topological space $M$.
A first candidate for equivariant cohomology might be
the singular cohomology of the orbit space $M/G$.
The example above of a circle $G=S^1$ acting on $M=S^2$ by rotation
shows that this is not a good candidate, since the orbit space $M/G$ is
a closed interval, a contractible space, so that its cohomology is
trivial.
In this example, we lose all information about the group action by passing to the quotient $M/G$.

A more serious deficiency of this example is that it is the quotient
of a nonfree action.
In general, a group action is said to be \emph{free} if the stabilizer of every point is the trivial subgroup.
It is well known that the orbit space
of a nonfree action is often ``not nice''---not
smooth or not Hausdorff.
However, the topologist has a way of turning every action into a free action without changing the homotopy type of the space.  
The idea is to find a contractible space $EG$ on which the group $G$ acts freely. 
Then $EG \times M$ will have the same homotopy type as $M$, and
no matter how $G$ acts on $M$, the diagonal action of $G$ on $EG
\times M$ will always be free.
The \emph{homotopy quotient} $M_G$ of $M$ by $G$,
also called the Borel construction, is defined to be the
quotient of $EG \times M$ by the diagonal action of $G$, and the
\emph{equivariant cohomology} $H_G^*(M)$ of $M$ is defined to be the
cohomology $H^*(M_G)$ of the homotopy quotient $M_G$.
Here $H^*(\ )$ denotes singular cohomology with any coefficient ring.

A contractible space on which a topological group $G$ acts freely is familiar from homotopy theory as the total space of a \emph{universal principal $G$-bundle} $\pi\colon EG \to BG$, of which every principal $G$-bundle is a pullback.
More precisely, if $P \to M$ is any principal $G$-bundle, then there is a map $f\colon M \to BG$, unique up to homotopy and called a \emph{classifying map} of $P \to M$, such that the bundle $P$ is isomorphic to the pullback bundle $f^*(EG)$.
The base space $BG$ of a universal bundle, uniquely defined up to homotopy equivalence, is called
the \emph{classifying space} of the group $G$.
The classifying space $BG$ plays a key role in equivariant cohomology, because it is the homotopy quotient of a point:
\[
\pt_G = (EG \times \pt)/G = EG/G = BG,
\]
so that the equivariant cohomology $H_G^*(\pt)$ of a point is the
ordinary cohomology $H^*(BG)$ of the classifying space $BG$.

It is instructive to see a universal bundle for the circle group.
Let $S^{2n+1}$ be the unit sphere in $\C^{n+1}$.
The circle $S^1$ acts on $\C^{n+1}$ by scalar multiplication.
This action induces a free action of $S^1$ on $S^{2n+1}$, and the quotient space is by definition the complex projective space $\C P^n$.
Let $S^{\infty}$ be the union $\bigcup_{n=0}^{\infty} S^{2n+1}$,
and let $\C P^{\infty}$ be the union 
$\bigcup_{n=0}^{\infty} \C P^n$. 
Since the actions of the circle on the spheres
are compatible with the inclusion of one sphere inside the next, there
is an induced action of $S^1$ on $S^{\infty}$. 
This action is free with quotient space $\C P^{\infty}$.
It is easy to see that all homotopy groups of $S^{\infty}$ vanish, for 
if a sphere $S^k$ maps into the infinite sphere $S^{\infty}$, then by
compactness its image lies in a finite-dimensional sphere $S^{2n+1}$.
If $n$ is large enough, any map from $S^k$ to $S^{2n+1}$ will be
null-homotopic.  
Since $S^{\infty}$ is a CW complex, the vanishing of all homotopy
groups implies that it is contractible.
Thus, the projection $S^{\infty} \to \C P^{\infty}$ is a universal
$S^1$-bundle and, up to homotopy equivalence, $\C P^{\infty}$ is the
classifying space $BS^1$ of the circle.

If $H^*(\ )$ is a cohomology functor, the constant map $M \to \pt$ from any space $M$ to a point induces a ring homomorphism $H^*(\pt) \to H^*(M)$,
which gives $H^*(M)$ the structure of a module over the  ring $H^*(\pt)$. 
Thus, the cohomology of a point serves as the coefficient ring in any cohomology theory.
For the equivariant real singular cohomology of a circle action, the coefficient ring is
\begin{align*}
H_{S^1}^*(\pt;\R) &= H^*\left(\pt_{S^1}; \R\right) = H^*(BS^1; \R) \\
&= H^*(\C P^{\infty};\R) \simeq \R [u],
\end{align*}
the polynomial ring generated by an element $u$ of degree $2$.
For the action of a torus $T= S^1 \times \cdots \times S^1 =
(S^1)^{\ell}$, the coefficient ring is the polynomial ring 
$H_T^*(\pt;\R) = \R[u_1,
\ldots, u_{\ell}]$, where each $u_i$ has degree $2$.

\section*{Equivariant Integration}

Let $G$ be a compact, connected Lie group.
Over a compact, oriented $G$-manifold, equivariant forms
can be integrated, but the values are in the coefficient ring
$H_G^*(\pt;\R)$, which is generally a ring of polynomials. 
According to Cartan, in the case of a circle action
 on a compact, oriented manifold,
an equivariant form of degree $2n$ is a sum 
\begin{multline} \label{e:combination}
\omega = \omega_{2n} + \omega_{2n-2}\, u 
+ \omega_{2n-4}\,u^2 \\
 + \cdots +
\omega_0\, u^n,
\end{multline}
where $\omega_{2j}\in \Omega^{2j}(M)^{S^1}$ is an $S^1$-invariant
$2j$-form on $M$.
If $\omega$ is closed under the Cartan differential, then it is called
an \emph{equivariantly closed extension} of the ordinary differential
form $\omega_{2n}$.
The equivariant integral $\int_M \omega$ is obtained by integrating
each $\omega_{2j}$ over $M$.
If $M$ has dimension $2n$, then the integral $\int_M \omega_{2j}$
vanishes except when $j=n$, and one has
\begin{align*}
\int_M \omega & = \int_M \omega_{2n} + \left(\int_M \omega_{2n-2}
\right) u \\
& \qquad \qquad + \cdots + \left( \int_M \omega_0 \right) u^n \\
& = \int_M \omega_{2n} + 0 + \cdots + 0 = \int_M \omega_{2n}.
\end{align*}

One peculiarity of equivariant integration is the possibility of
obtaining a nonzero answer while integrating a form over a manifold
whose dimension is not equal to the degree of the form.
For example, if $M$ has dimension $2n-2$ instead of $2n$, then the
integral over $M$ of the equivariant $2n$-form $\omega$ above is
\[
\int_M \omega = \left( \int_M \omega_{2n-2} \right) u,
\]
since for dimensional reasons all other terms are zero.
From this, one sees that an equivariant integral for a circle action
is in general not a real number, but a polynomial in $u$.

\section*{Localization}

What kind of information can be mined from the fixed points of an
action?
If a Lie group $G$ acts smoothly on a manifold,
then for each $g \in G$, the action induces a diffeomorphism
$\ell_g\colon M \to M$.
At a fixed point $p \in M$, the differential $\ell_{g*}\colon T_pM \to
T_pM$ is a linear automorphism of the tangent space,
giving rise to a representation of the group $G$ on the tangent space
$T_pM$. 
Invariants of the representation are then invariants of the action at
the fixed point.
For a circle action, at an isolated fixed point $p$,
the tangent space $T_pM$ decomposes into a direct sum
$L^{m_1} \oplus \cdots \oplus L^{m_n}$, where $L$ is the
standard representation of the circle on the complex plane $\C$
and $m_1, \ldots, m_n$ are nonzero integers.
The integers $m_1, \ldots, m_n$ are called the \emph{exponents}
of the circle action at the fixed point $p$.
They are defined only up to sign, but if $M$ is oriented,
the sign of the product $m_1 \cdots m_n$ is well defined by the
orientation of $M$.

When a torus $T=(S^1)^{\ell}$ acts on a compact, oriented manifold $M$ with isolated fixed point set $F$,
for any closed $T$-equivariant form $\omega$ on $M$,
the equivariant integration formula states that
\begin{equation}\label{e:localization}
\int_M \omega = \sum_{p \in F} \frac{i_p^*\omega}{e^T(\nu_p)},
\end{equation}
where $i_p^*\omega$ is the restriction of the equivariantly closed
form $\omega$ to a fixed 
point $p$ and $e^T(\nu_p)$ is the equivariant Euler class of the
normal bundle $\nu_p$ to $p$ in $M$.
Of course, the normal bundle to a point $p$ in a manifold
$M$ is simply the tangent
space $T_pM$, but formula~\eqref{e:localization} is
stated in a way to allow for easy generalization:
when $F$ has positive-dimensional components,
the sum over the fixed points is replaced by an
integral over the components $C$ of the fixed point set
and $\nu_p$ is replaced by $\nu_C$, the normal bundle
to the component $C$.
In formula~\eqref{e:localization},
the degree of the form $\omega$ is not assumed to be
equal to the dimension of the manifold $M$, and so the left-hand side
is a polynomial in $u_1, \ldots, u_{\ell}$, while the right-hand side
is a sum of rational expressions in $u_1, \ldots, u_{\ell}$,
and it is part of the theorem that the equivariant Euler classes
$e^T(\nu_p)$ are nonzero and that there will be cancellation on the
right-hand side so that the sum becomes a polynomial.

Return now to a circle action with isolated fixed points on a compact,
oriented manifold $M$ of dimension $2n$.
Let $\omega$ be a closed equivariant form of degree $2n$ on $M$.
Since the restriction of a form of positive degree to a point is zero,
on the right-hand side of \eqref{e:localization} all terms in $\omega$
except
$\omega_0\, u^n$ restrict to zero at a fixed point $p\in M$:
\begin{align*}
i_p^* \omega &= \sum_{j=0}^n (i_p^* \omega_{2n-2j})\, u^j \\
& = (i_p^* \omega_0)\, u^n = \omega_0(p)\, u^n.
\end{align*}
The equivariant Euler class
$e^{S^1}(\nu_p)$ turns out to be $m_1 \cdots m_n\, u^n$,
where $m_1, \ldots, m_n$ are the exponents of the circle
action at the fixed point $p$.
Therefore, the equivariant integration formula for a circle action
assumes the form
\[
\int_M \omega_{2n} = \int_M \omega = \sum_{p \in F}
\frac{\omega_0}{m_1 \cdots m_n} (p).
\]
In this formula, $\omega_{2n}$ is an ordinary differential form of
degree $2n$ on $M$, $\omega$ is an equivariantly closed extension of
$\omega_{2n}$, and $\omega_0$ is the coefficient of the $u^n$ term in
$\omega$ as in \eqref{e:combination}. 

\section*{Applications}

In general, an integral of an ordinary differential form on a
compact, oriented manifold can be computed as a finite sum using the
equivariant integration formula if the manifold has a torus action
with isolated fixed points and the form has an equivariantly closed
extension. 
These conditions are not as restrictive as they seem,
since many problems come naturally with the action of a compact Lie
group, and one can always restrict the action to that of a maximal
torus.
It makes sense to restrict to a maximal torus, instead of
any torus in the group, because the larger the torus, the smaller the
fixed point set, and hence the easier the computation.

As for the question of whether a form has an equivariantly closed
extension, in fact a large collection of forms automatically do.  
These include characteristic classes of vector bundles on a
manifold.
If a vector bundle has a group action compatible with the group action
on the manifold, then the equivariant characteristic classes of the
vector bundle will be equivariantly closed extensions of its ordinary
characteristic classes.

A manifold on which every closed form has an equivariantly closed extension is said to be \emph{equivariantly formal}.  Equivariantly formal manifolds include all manifolds whose cohomology vanishes in odd degrees.  In particular, a homogeneous space $G/H$, where $G$ is a compact Lie group and $H$ is a closed subgroup of maximal rank, is
equivariantly formal.

The equivariant integration formula is a powerful tool for computing
integrals on a manifold.
If a geometric problem with an underlying torus action can be
formulated in terms of integrals, then there is a good chance that the
formula applies.
For example, it has been applied to show that the stationary phase
approximation formula is exact for a symplectic action (Atiyah--Bott 1984),
to calculate the number of rational curves in a quintic threefold
(Kontsevich 1995, Ellingsrud--Str\o mme 1996),
to calculate the characteristic numbers of a compact homogeneous
space (Tu 2010), and to derive the Gysin formula for a fiber bundle with
homogenous space fibers (Tu preprint 2011).
In the special case where the vector field $X$ is generated by
a circle action,
the Gauss--Bonnet--Hopf formula~\eqref{e:gbh} is
a consequence of the equivariant integration formula.
Equivariant cohomology has also helped to elucidate
the work of Witten on supersymmetry, Morse 
theory, and Hamiltonian actions (Atiyah--Bott 1984, Jeffrey--Kirwan 1995).

The formalism of equivariant cohomology carries over from singular
cohomology to other cohomology theories such as $K$-theory,
Chow rings, and quantum cohomology.
There are similar localization formulas that compare
the equivariant functor of a $G$-space to that of the fixed
point set of $G$ or of some subgroup of $G$ (for example,
Segal 1968 and Atiyah--Segal 1968).
In the fifty years since its inception, equivariant cohomology has
found applications in topology, differential geometry,
symplectic geometry, algebraic geometry, $K$-theory,
representation theory, and combinatorics, among other fields,
and is currently a vibrant area of research.

\section*{Acknowledgments}
This article is based on a talk given at the National Center for Theoretical Sciences, National Tsing Hua University, Taiwan.
The author gratefully acknowledges helpful discussions with Alberto Arabia, Aaron W.~Brown, Jeffrey D.~Carlson, George Leger, and Winnie Li during the preparation of this article, as well as the support of the Tufts Faculty Research Award Committee in 2007--2008.

\renewcommand{\refname}{For Further Reading}

\end{multicols}
\end{document}